\documentclass[12pt]{article}
\usepackage[utf8]{inputenc}
\usepackage{cite}
\usepackage{graphicx}
\usepackage{amsmath}
\usepackage{amsthm}
 \usepackage{amssymb}
 \usepackage{geometry}
 \usepackage[nottoc]{tocbibind}
 \usepackage{listings}
\usepackage{pdfpages}
\usepackage{algorithm}
\usepackage{booktabs}
\usepackage{algpseudocode}
\usepackage{mathtools}
\usepackage{hyperref}
\theoremstyle{definition}

\newtheorem{claim}{Claim}[section]

\theoremstyle{remark}

\usepackage{pdfpages}
\providecommand{\keywords}[1]
{

  \textbf{Key words. }#1
}

\title{A Monotone, Second Order Accurate Scheme\\ for Curvature Motion}
\author{Selim Esedo\=glu\thanks{Email: esedoglu@umich.edu} \\ University of Michigan \and Jiajia Guo\thanks{Email: jiajiag@umich.edu (Corresponding author)}\\ University of Michigan}
\date{\today}

\begin{document}

\maketitle
\begin{abstract}
\noindent We present a second order accurate in time numerical scheme for curve shortening flow in the plane that is unconditionally monotone.
It is a variant of threshold dynamics, a class of algorithms in the spirit of the level set method that represent interfaces implicitly.
The novelty is monotonicity: it is possible to preserve the comparison principle of the exact evolution while achieving second order in time consistency.
As a consequence of monotonicity, convergence to the viscosity solution of curve shortening is ensured by existing theory.
\end{abstract}


\section{Introduction}
In this short note, we report a second order accurate threshold dynamics algorithm for simulating curvature motion (curve shortening) in the plane that is {\em monotone}: It respects the comparison principle of the exact evolution.
Existing theory \cite{barles_souganidis, ishii_pires_souganidis} then immediately implies that the approximate evolution generated by the scheme converges to the viscosity solution of mean curvature motion under appropriate conditions.

The finding is surprising, as previous studies e.g. \cite{ruuththesis, heintz05, zaitzeff1} that explored the idea of designing high order accurate versions of threshold dynamics speculated that monotonicity may need to be sacrificed.
This note shows that, at least in two dimensions, this need not be so.
To our knowledge, the algorithm presented is the only one of its kind (level-set style numerical algorithm capable of handling topological changes implicitly) that is rigorously shown to be unconditionally monotone, and consistent to second order, and thereby convergent.
All the advantages of the original scheme are retained, as the version proposed here differs only in its choice of convolution kernel, replacing the standard choice of Gaussian with a carefully chosen linear combination of Gaussians.
Hence, at least in two dimensions, there is a very special choice of a convolution kernel.

\section{Previous work}
There are several relevant contributions to high order in time versions of threshold dynamics in existing literature.
The first contribution is from the Ph.D. thesis of Ruuth \cite{ruuththesis}.
There, a second order accurate, multistep scheme inspired by Richardson extrapolation is proposed that is numerically demonstrated to achieve second order accuracy in time in two and three dimensional examples, at least while the evolving interface remains smooth.
However, the stability of that scheme (whether by maximum principles or energy methods) appears hard to study, and there are no results to that effect (or even a careful consistency calculation) available.

The topic of high order accurate schemes for curvature motion also comes up naturally as a byproduct in studies focused on adapting threshold dynamics to high order geometric motions such as Willmore flow \cite{grzibovskis_heintz, esedoglu_ruuth_tsai2}, where the idea of using linear combinations of Gaussians as the convolution kernel to cancel out undesirable terms in consistency calculations plays the same prominent role.
In \cite{heintz05}, designing a second order accurate in time version of threshold dynamics for curvature motion by a judicious choice of convolution kernel is floated, but the proposed approach would merely result in a second order accurate evaluation of curvature of the interface at the beginning of a time step, which is different from (and not sufficient for) advancing the interface under curvature motion with the requisite third order local truncation error.
In reality, such a scheme would still be just first order consistent.
Moreover, it is stated that the resulting kernel would not be positive everywhere, and thus the resulting scheme would violate monotonicity.

More recently, second order versions of threshold dynamics are proposed in \cite{zaitzeff1}.
One is multi-step, similar to that of Ruuth in \cite{ruuththesis}, and therefore not likely to be monotone.
However, unlike in \cite{ruuththesis}, it comes with a careful consistency calculation, which verifies second order consistency (in addition to numerical evidence) in two and three dimensions.
The other proposed scheme of \cite{zaitzeff1} is multi-stage, and therefore also unlikely to be monotone.
It is, however, second order consistent in two dimensions, and most notably, satisfies an unconditional energy stability property in any dimension: it dissipates the Lyapunov functional for threshold dynamics discovered in \cite{esedoglu_otto}.
In the broader context of level-set type methods that represent interfaces implicitly, the early contribution \cite{walkington} proposes second order, energy (total variation) diminishing schemes for the level set formulation of mean curvature motion, but reports difficulties with (slow or lack of convergence of iterative solvers on) the nonlinear algebraic systems that need to be solved at every time step.

It is also worth recalling that using different (namely, in this case, non-radially symmetric) kernels in threshold dynamics comes up in its extensions to anisotropic curvature flows \cite{ruuth_merriman, bbc, elsey_esedoglu_anisotropy, ejz}.
In particular, barrier type theorems \cite{elsey_esedoglu_anisotropy, ejz} show that any threshold dynamics scheme that is at all consistent (never mind second order) with certain anisotropic curvature flows in three dimensions cannot possibly be monotone.

\section{The standard algorithm}
Recall that threshold dynamics algorithm of Merriman, Bence, and Osher \cite{mbo92, merriman_bence_osher} generates a discrete in time approximation to the motion by mean curvature of an interface $\partial\Sigma^0$ given as the boundary of an initial set $\Sigma^0 \subset \mathbb{R}^d$ by alternating the two steps of convolution and thresholding:
\begin{center}
\fbox{\begin{minipage}[t]{0.8\columnwidth}
\textbf{Algorithm:}(MBO'92): Given a time step size $t>0$, alternate the following steps:
\begin{enumerate}
\item Convolution:
\begin{equation}
\psi^{k} = K_t * \mathbf{1}_{\Sigma^{k}}\label{eq:mbo1}
\end{equation}
\item Thresholding:
\begin{equation}
\label{eq:mbo2}
\Sigma^{k+1}=\left\{ x\,:\,\psi^{k}(x) \geq \lambda \right\} .
\end{equation}
\end{enumerate}
\end{minipage}}
\end{center}
where we write
\begin{equation*}
K_t(x) = \frac{1}{t^{\frac{d}{2}}} K \left( \frac{x}{\sqrt{t}} \right)
\end{equation*}
for a smooth convolution kernel $K:\mathbb{R}^d\to\mathbb{R}$ of total mass $2\lambda>0$ and sufficiently rapid decay at $|x| \to \infty$.
The kernel $K$ was chosen in \cite{mbo92} originally to be the Gaussian:
\begin{equation}
\label{eq:gaussian}
G(x) = \frac{1}{(4\pi)^\frac{d}{2}} \exp\left(-\frac{|x|^2}{4}\right)
\end{equation}
but choosing it something else was also raised as a possibility in the same work.
With choice (\ref{eq:gaussian}), convergence of scheme (\ref{eq:mbo1}) \& (\ref{eq:mbo2}) had been established in a number of previous studies, including \cite{evans_mbo, barles_georgelin, ishii_pires_souganidis, swartz}.
There are even convergence results \cite{laux_otto, laux_otto2} in the multiphase setting.

In this note, we restrict attention to radially symmetric convolution kernels of the form
\begin{equation*}
K (x) = \sum_{j=1}^N c_j G_{\alpha_j} (x)
\end{equation*}
and ask whether the coefficients $\alpha_j$ and $c_j$ can be chosen so that
\begin{enumerate}
    \item $K(x)\geq 0$, and
    \item Scheme (\ref{eq:mbo1}) \& (\ref{eq:mbo2}) is second order consistent.
\end{enumerate}
We are surprised to find out that the answer is yes when $d=2$.
Recall that $K(x)\geq 0$ implies unconditional monotonicity of the scheme:
\begin{equation*}
\Omega^0 \subset \Sigma^0 \Longrightarrow \Omega^k \subset \Sigma^k \mbox{ for all } k=1,2,3,\ldots
\end{equation*}
regardless of the time step size $k>0$.
Preserving this fundamental qualitative feature of the exact evolution is tremendously helpful in establishing stability and convergence of numerical schemes.
\section{A special kernel in dimension $d=2$}
\label{sec:2d}
In this section, we carefully exhibit a positive convolution kernel that endows scheme (\ref{eq:mbo1}) \& (\ref{eq:mbo2}) with second order consistency and, thanks to positivity, monotonicity.
Assume that the initial interface is given as the graph of a smooth function $g:\mathbb{R} \to \mathbb{R}$ with $g(0)=0$ and $g'(0)=0$, and the initial set is $\Sigma^0 = \{ (x,y) \, : \, y \geq g(x) \}$.
The exact solution $y = \phi(x,t)$ of curvature motion solves the PDE
\begin{equation} \label{pde}
\left\{
    \begin{split}\textbf{}
      &  \phi_t = \frac{\phi_{xx}}{1+\phi_x^2}\\
      &  \phi(0,x) = g(x)   
    \end{split}
\right.
\end{equation}
Taylor expanding $\phi(0,t)$ in $t$ at $t=0$ and converting all time derivatives to spatial ones via the equation, one gets:
\begin{equation}
\label{eq:exact2d}
\phi(0,t) = t \, g''(0) + t^2 \left( \frac{1}{2} g^{(iv)}(0) - \big( g''(0) \big)^3(0) \right) + O(t^3)
\end{equation}
as $t\to 0^+$.
We will demand that one step of the threshold dynamics scheme gives an interface that crosses the $y$-axis at $\phi(0,t)$ up to $O(t^3)$ terms.

Taylor expansion for the interface after one step of threshold dynamics using a Gaussian kernel had been obtained in multiple previous studies, e.g. \cite{mascarenhas, ruuththesis, ruuth0, grzibovskis_heintz}.
The first step in those calculations is to expand the convolution step (\ref{eq:mbo1}) of the algorithm along the $y$-axis, which is of course linear in the kernel $K$.
Define
\begin{equation}
\label{eq:theta}
\theta(p) = \sum_{j=1}^N \alpha_j^{\frac{p}{2}} c_j.
\end{equation}
We get
\begin{equation}
\label{eq:expansion1}
\begin{split}
K * \mathbf{1}_\Sigma(0,y) =& \frac{1}{2} \theta(0) - \frac{1}{2\sqrt{\pi}} \frac{y}{\sqrt{t}} \theta(-1) + \frac{1}{24\sqrt{\pi}} \frac{y^3}{t^{3/2}} \theta(-3)\\
&+ \frac{g''(0)}{2\sqrt{\pi}} \sqrt{t} \theta(1) - \frac{g''(0)}{8\sqrt{\pi}} \frac{y^2}{\sqrt{t}} \theta(-1)\\
&+ \frac{g^{(iv)}(0)}{4\sqrt{\pi}} t^{3/2} \theta(3) + \frac{3 \big( g''(0) \big)^2}{8\sqrt{\pi}}  \sqrt{t} y \theta(1)\\
&- \frac{5 \big( g''(0) \big)^3}{8\sqrt{\pi}} t^{3/2} \theta(3) + O(t^{5/2})
\end{split}
\end{equation}
under the assumption that $y=O(t)$ as $t\to 0$.
The next step is to obtain the expansion for the interface after the thresholding step (\ref{eq:mbo2}); at that point, the dependence of the scheme on the convolution kernel $K$ is no longer linear.
We substitute the ansatz $y = a_1 t + a_2 t^2 + O(t^3)$ into (\ref{eq:expansion1}) to get
\begin{equation}
\label{eq:Kansatz}
K_t * \mathbf{1}_{\Sigma}(0,y) = A_0 + A_1 \sqrt{t} + A_3 t^{3/2} + O(t^{5/2})
\end{equation}
where
\begin{equation}
\label{eq:A0}
A_0 = \frac{1}{2} \theta(0),
\end{equation}
and
\begin{equation}
\label{eq:A1}
A_1 = \frac{1}{2\sqrt{\pi}} 
\Big( 
\theta(1) g''(0)
-a_1 \theta(-1)
\Big),
\end{equation}
and
\begin{equation}
\label{eq:A3}
\begin{split}
A_3 =& \frac{1}{8\sqrt{\pi}} \theta(3) \left( 2 g^{(iv)}(0) - 5 \big( g''(0) \big)^3 \right)
+ \frac{3 \big(g''(0)\big)^2}{8\sqrt{\pi}} \theta(1) a_1\\
-& \frac{g''(0)}{8\sqrt{\pi}} \theta(-1) a_1^2
+ \frac{1}{24\sqrt{\pi}} \theta(-3) a_1^3
- \frac{1}{2\sqrt{\pi}} \theta(-1) a_2
\end{split}
\end{equation}
It turns out that taking $N=3$ is sufficient for our purposes in this section.
Thus, for the rest of this section, we take our kernel $K$ to be of the form
\begin{equation}
\label{eq:K2d}
K = G + c_1 G_{\alpha_1} + c_2 G_{\alpha_2}
\end{equation}
From (\ref{eq:A0}), we see that the convolution level will be given by
\begin{equation}
\label{eq:lambda}
\lambda = \frac{1+c_1+c_2}{2}.
\end{equation}
Setting (\ref{eq:Kansatz}) equal to (\ref{eq:lambda}), we require $A_1=0$ and $A_2=0$.
At the $O(\sqrt{t})$ level, solving $A_1=0$ for $a_1$, we find:
\begin{equation}
\label{eq:a1}
a_1 = \frac{\theta(1)}{\theta(-1)} g''(0).
\end{equation}
Having determined $a_1$, $A_3$ can now be expressed as
\begin{equation}
\begin{split}
A_3 =& \frac{\theta(3)}{4\sqrt{\pi}} g^{(iv)}(0) + \frac{\theta^3(1)\theta(-3) + 6 \theta^2(1)\theta^2(-1) - 15 \theta^3(-1)\theta(3)}{24\sqrt{\pi}\theta^3(-1)} \big(g''(0)\big)^3\\
&- \frac{\theta(-1)}{2\sqrt{\pi}} a_2
\end{split}
\end{equation}
Setting $A_3=0$ and solving for $a_2$, we find
\begin{equation}
\label{eq:a2}
a_2 = \frac{\theta(3)}{2\theta(-1)} g^{(iv)}(0) + \frac{\theta^3(1)\theta(-3) + 6 \theta^2(1)\theta^2(-1) - 15 \theta^3(-1)\theta(3)}{12\theta^4(-1)} \big(g''(0)\big)^3.
\end{equation}
We can now compare (\ref{eq:a1}) \& (\ref{eq:a2}) with (\ref{eq:exact2d}).
To match the two expansions for some effective choice of time step size in Algorithm (\ref{eq:mbo1}) \& (\ref{eq:mbo2}), we need:
\begin{equation}
\label{eq:theta1}
\frac{\theta^2(1)}{\theta^2(-1)} = \frac{\theta(3)}{\theta(-1)}
\end{equation}
and
\begin{equation}
\label{eq:theta2}
\frac{\theta(3)}{\theta(-1)} = - \frac{\theta^3(1)\theta(-3) + 6 \theta^2(1)\theta^2(-1) - 15 \theta^3(-1)\theta(3)}{12\theta^4(-1)} 
\end{equation}
along with the proviso $\theta(-1)\not= 0$ that we will verify at the end.
Taking $\alpha_1 = 4$ and $\alpha_2 = 1/4$, equation (\ref{eq:theta1}) becomes
\begin{equation}
72 c_1 + 18 c_2 + 225 c_1 c_2 = 0.
\end{equation}
which gives
\begin{equation}
\label{eq:c2}
c_2 = -\frac{8 c_1}{25 c_1 + 2}.
\end{equation}
Substituting into (\ref{eq:theta2}), we get
\begin{equation}
\label{eq:c1}
\frac{(10c_1+1)^2 \big(1000 c_1^3 - 2175 c_1^2 + 210 c_1 + 64\big)}{12 (5c_1+2)^5} = 0
\end{equation}
The polynomial
\begin{equation}
\label{eq:polynomial}
p(x) = 1000 c_1^3 - 2175 c_1^2 + 210 c_1 + 64
\end{equation}
satisfies $p(1/5) = 27$ and $p(1/4) = -61/16$, and hence has a root in $(\frac{1}{5},\frac{1}{4})$, at which the denominator of (\ref{eq:c1}) does not vanish.
Taking this root as the value of $c_1$, i.e.
\begin{equation}
\label{eq:c1value}
c_1 \approx 0.2444098,
\end{equation}
equation (\ref{eq:c1}) is then satisfied.
Substituting into (\ref{eq:c2}) determines $c_2$:
\begin{equation}
\label{eq:c2value}
c_2 \approx -0.2410874.
\end{equation}
We note that $\theta(-1) = 1 + \frac{1}{2}c_1 - 2c_2 \not= 0$, as hoped for.

\begin{figure}[h]
\begin{center}
\includegraphics[scale=0.75]{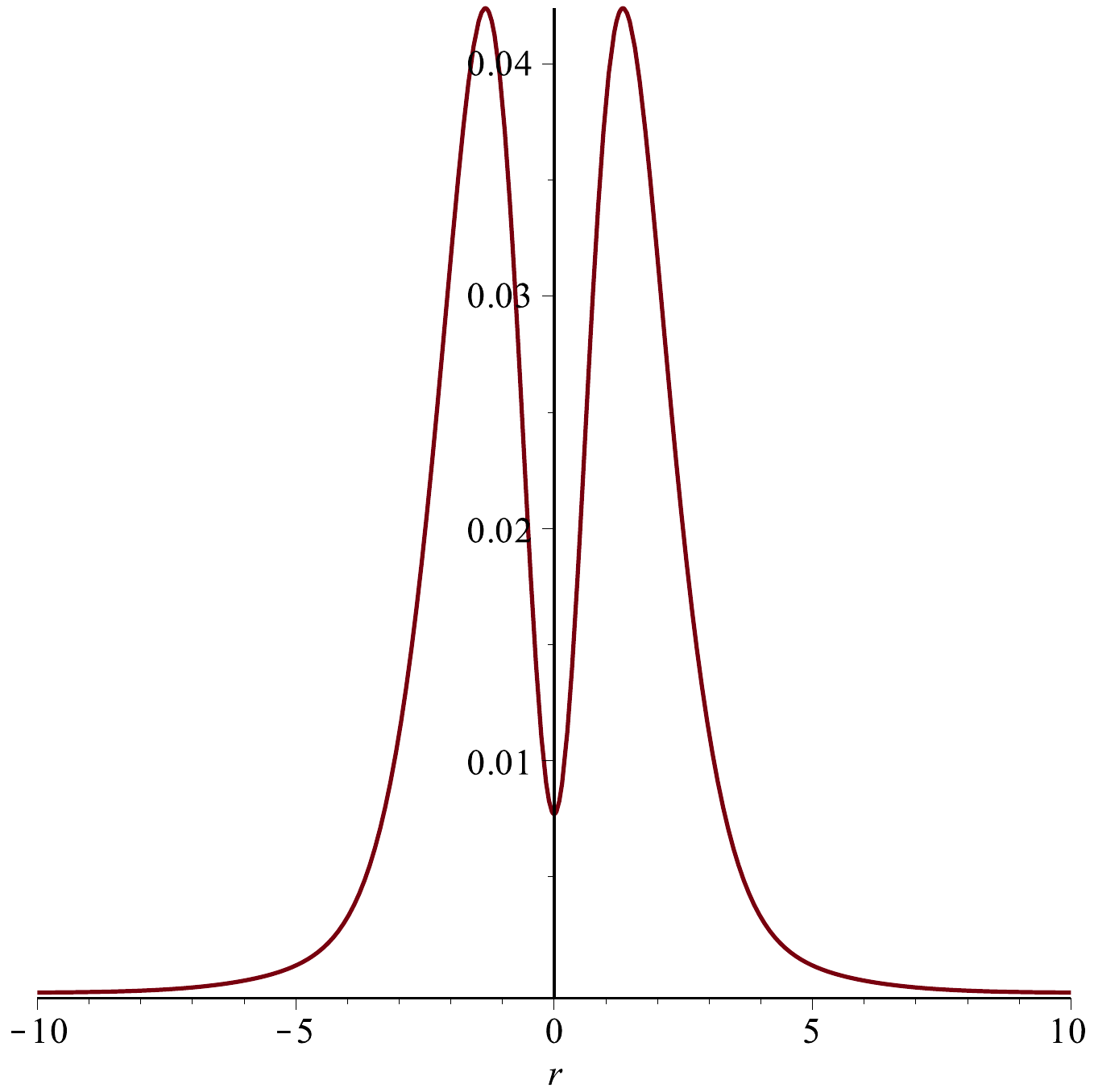}
\caption{\footnotesize Radial profile of the special convolution kernel $K$ that endows the standard threshold dynamics algorithm of Merriman, Bence, and Osher with second order accuracy in two dimensions, as well as with monotonicity.}
\label{fig:kernel}
\end{center}
\end{figure}

Returning to (\ref{eq:a1}), we see that when the convolution kernel $K$ in Step (\ref{eq:mbo1}) is given by (\ref{eq:K2d}) with $\alpha_1=4$, $\alpha_2 = \frac{1}{4}$ and the two coefficients $c_1, c_2$ are given by (\ref{eq:c1value}) \& (\ref{eq:c2value}), Algorithm (\ref{eq:mbo1}) \& (\ref{eq:mbo2}) is {\em second order accurate} with the (rescaled) effective time step size
\begin{equation}
\tau = \frac{\theta(1)}{\theta(-1)} t = \frac{1+2c_1+\frac{1}{2}c_2}{1+\frac{1}{2}c_1+2c_2} t \approx 2.137831 t.
\end{equation}
Figure \ref{fig:kernel} shows a plot of the radial profile of the kernel.
It appears to be positive; we now show that it indeed is.
Let $\xi = \exp(-\frac{1}{16}r^2)$ so that $\xi\in(0,1]$.
Then,
$$ K = \frac{1}{\pi} \xi q(\xi) \mbox{ where } q(\xi) = c_2 \xi^{15} + \frac{1}{4} \xi^3 + \frac{1}{16} c_1.$$
We have
\begin{equation}
\begin{split}
q(\xi) &= \xi^3 \left( \frac{1}{4} + c_2 \xi^{12} \right) + \frac{1}{16} c_1\\
&\geq \xi^3 \left( \frac{1}{4} + c_2 \right) + \frac{1}{16} c_1\\
&\geq 0
\end{split}
\end{equation}
since $c_1\in(\frac{1}{5},\frac{1}{4})$ and so $c_2\in(-\frac{8}{33},-\frac{8}{35})$.
We have established the following:
\begin{claim}
Let $c_1$ be the root of the polynomial (\ref{eq:polynomial}) in $(\frac{1}{5},\frac{1}{4})$.
Let $c_2$ be given in terms of $c_1$ by (\ref{eq:c2}).
Let $K$ be the convolution kernel
$$ K = G + c_2 G_4 + c_3 G_{\frac{1}{4}}. $$
Then, scheme (\ref{eq:mbo1}) \& (\ref{eq:mbo2}) is monotone, and second order consistent with curvature motion in the plane.
The discrete in time evolutions generated by the scheme (extended from sets to functions in the natural way of e.g. \cite{ishii_pires_souganidis}) converge uniformly to the unique viscosity solution \cite{evans_spruck, chen_giga_goto} of curvature motion on any finite time interval, starting from bounded, uniformly continuous initial data.
\end{claim}

\begin{proof}
Immediate consequence of positivity, smoothness, and decay properties of the kernel, the consistency calculation above, and the theory of \cite{barles_souganidis, ishii_pires_souganidis}.
In fact, our kernel satisfies the conditions of \cite{ishii_pires_souganidis}.
\end{proof}

\section{Numerical demonstration}
We demonstrate that second order accuracy (in time) is indeed achieved by Algorithm (\ref{eq:mbo1}) \& (\ref{eq:mbo2}) using the convolution kernel (\ref{eq:K2d}).
To that end, and to minimize any potential issues with insufficient spatial resolution, we implement the algorithm {\bf 1.} in the radial case to test on a shrinking circle, and {\bf 2.} in case the interface is given as the periodic graph of a function.

\begin{figure}[h]
\centering
\includegraphics[width=6in]{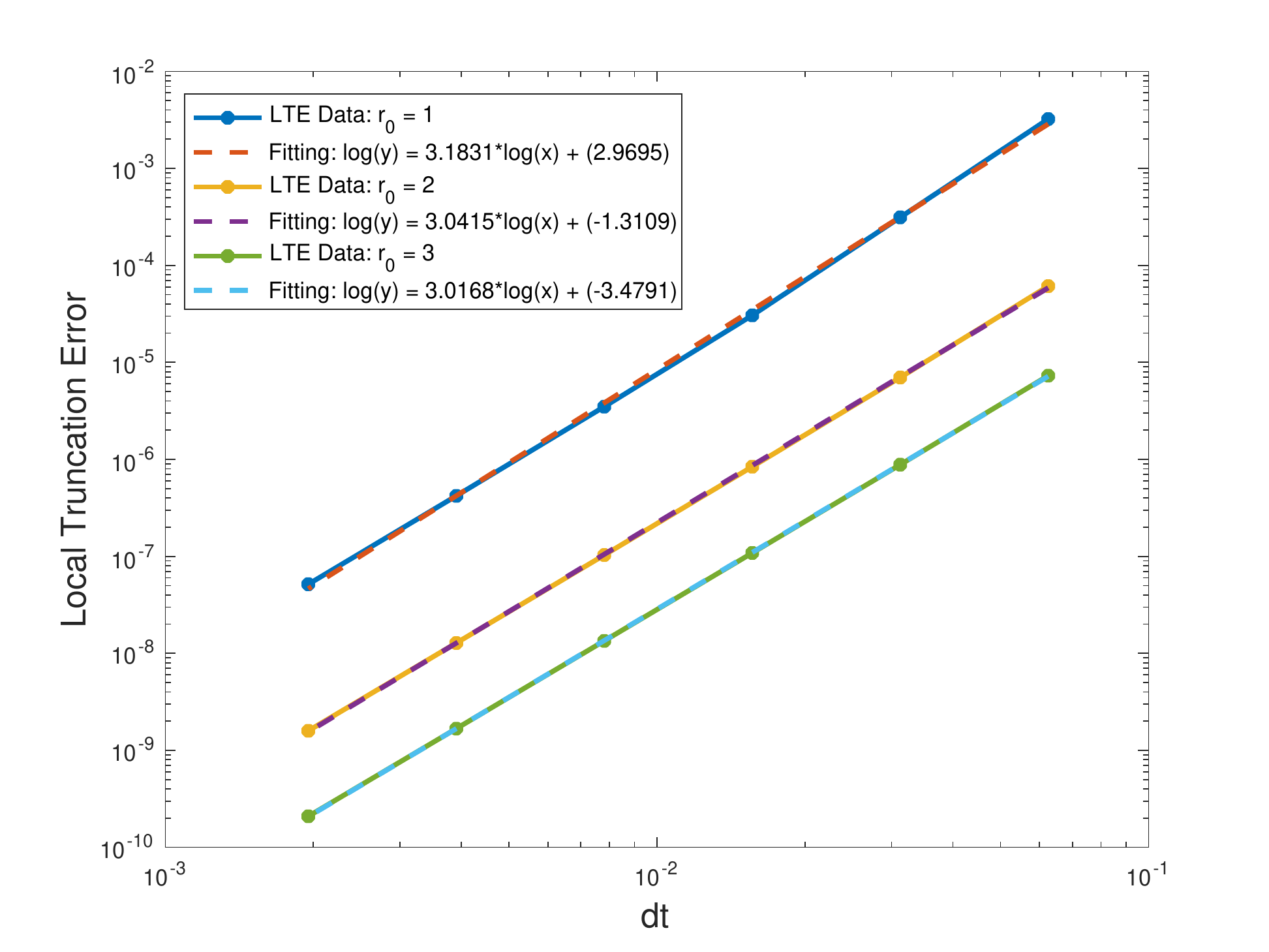}   
\caption{\footnotesize Local truncation error, for initial condition given by circles of radii $r_0=1,2,3$.}
\label{fig:LTE}
\end{figure}

For the shrinking circle test, we merely test the local truncation error by taking a single time step (of various sizes) with the algorithm, starting from initial radii of $r_0=1,2$, and $3$.
In this case, the exact solution of curvature motion is given by $\sqrt{r_0^2 - 2t}$.
The convolution in the algorithm is calculated very accurately in polar coordinates.
The expected $O(t^3)$ rate of decay of the local error is observed, as show in Figure \ref{fig:LTE}.

\begin{figure}[h] 
\centering
\includegraphics[width=6in]{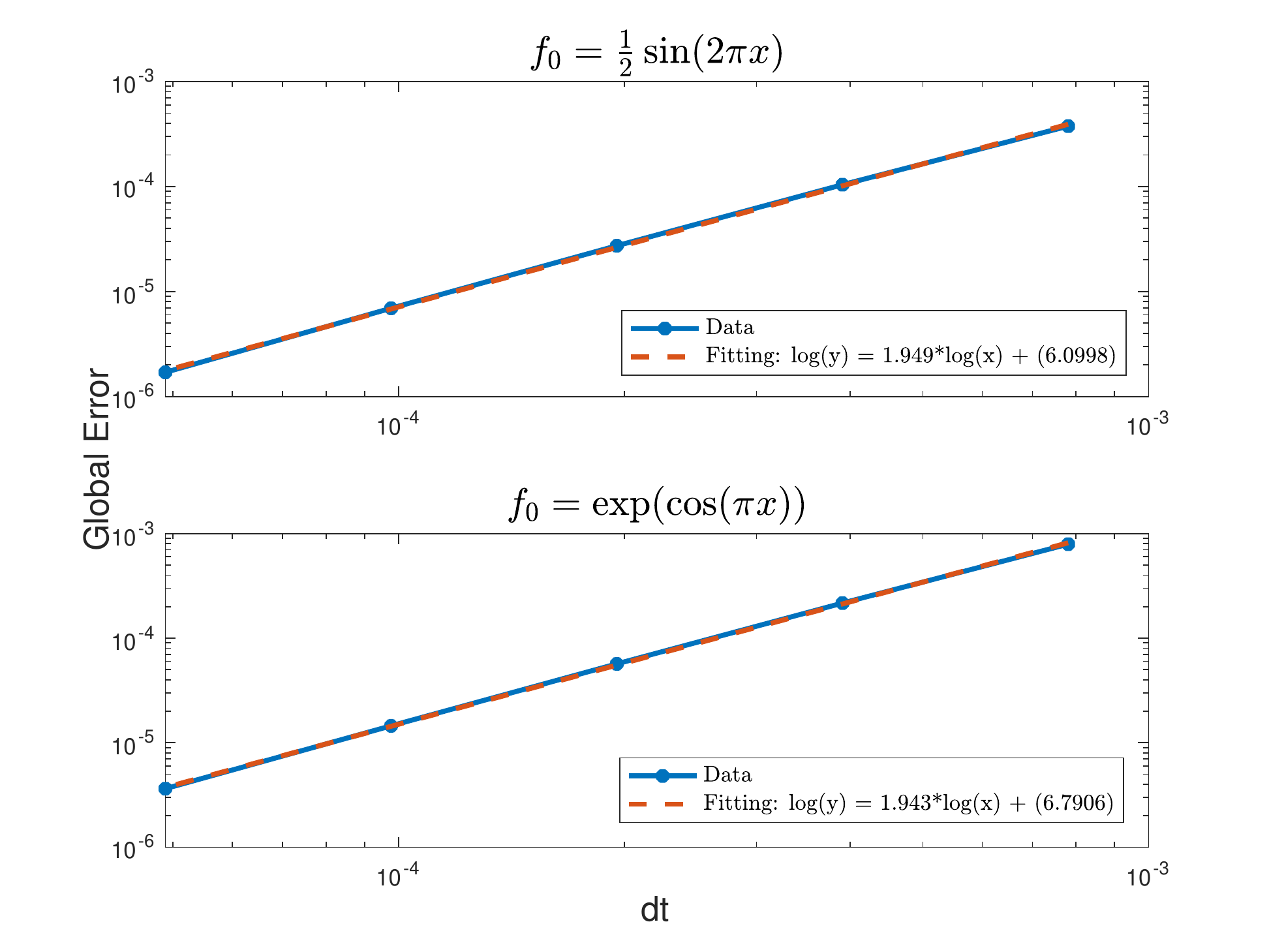}   
\caption{\footnotesize Global error, for initial condition given $f_0 = \frac{1}{2}\sin(2\pi x)$ and $f_0 = \exp\left(\cos(\pi x)\right)$.}
\label{fig:Globalerror}
\end{figure}

To test on interfaces given as graphs of functions, we measure the global error at final time $T=1/40$, starting from periodic initial conditions $y = f_0(x)$ where $f_0(x) = \frac{1}{2}\sin(2\pi x)$ and $f_0(x) = \exp(\cos(\pi x))$.
The benchmark solution of the PDE (\ref{pde}) is obtained by an extremely fine finite differences discretization (forward Euler time steps, and centered differences in space).
The algorithm is implemented by discretizing the $x$-axis, and for each discrete $x$-value, finding the value of $y$ for which the convolution integral
\begin{equation}
\label{eq:integral}
\begin{split} 
      G_{t} * \mathbf{1}_{\Sigma}(x,y) &= \int_{-\infty}^{+\infty}G_{t}(x - \bar{x}) \int_{-\infty}^{f_0(\bar{x})} \frac{1}{\sqrt{4\pi \delta t}} e^{-\frac{(y-\bar{y})^2}{4\delta t}}d\bar{y}d\bar{x}\\
      &=\frac{1}{2} + \frac{1}{2}\int_{-\infty}^{+\infty} G_{t}(x-\bar{x}) \, \mbox{erf}\left(\frac{f(\bar{x})-y}{2\sqrt{t}}\right)d\bar{x}
\end{split}
\end{equation}
equals the thresholding value $\lambda$, where the interface at the current time step is represented by $y=f(x)$.
The convolution integral (\ref{eq:integral}) is estimated numerically, truncating its integrand once it falls below a tolerance.
The expected $O(t^2)$ scaling of the global error can be seen in the Tables \ref{table:1} \& \ref{table:2} and Figure \ref{fig:Globalerror}. 
\begin{table}[h!]
\begin{center}
\begin{tabular}{c c c c c c} 
 \hline
 Number of time steps & 32 & 64 & 128 & 256 & 512\\
 \hline
 $L^2$ error & 3.76e-04 & 1.04e-04 & 2.73e-05  & 6.93e-06 & 1.70e-06\\ 
 \hline
 Order & - & 1.9 & 1.9  & 2.0 & 2.0
\end{tabular}
\end{center}
\caption{Error and order for $f_0(x) = \frac{1}{2}\sin(2\pi x)$.}
\label{table:1}
\end{table}
\begin{table}[h!]
\begin{center}
\begin{tabular}{c c c c c c} 
 \hline
 Number of time steps & 32 & 64 & 128 & 256 & 512\\
 \hline
 $L^2$ error & 7.92e-04 & 2.16e-04 & 5.67e-05  & 1.45e-05 & 3.64e-06\\ 
 \hline
 Order & - & 1.9 & 1.9  & 2.0 & 2.0
\end{tabular}
\end{center}
\caption{Error and order for $f_0(x) = \exp(\cos(\pi x))$.}
\label{table:2}
\end{table}

\section{Dimension $d=3$}
\label{sec:3d}
In this section, we show that the simple kernel construction of Section \ref{sec:2d} as a linear combination of Gaussians will not work in higher dimensions.
This was mentioned in \cite{zaitzeff1} in passing; here we give a careful explanation.
For $d=3$, assume that the initial interface $\partial\Sigma^0$ is given as the graph of a smooth function $g:\mathbb{R}^2 \to \mathbb{R}$ with $g(0,0) = 0$ and $\nabla g(0,0) = 0$ so that $\Sigma^0 = \{ (x,y,z) \, : \, z \geq g(x,y) \}$.
The exact solution of motion by mean curvature is described at least for short time by the PDE:
\begin{equation}
\label{eq:MCM3d}
\left\{
\begin{split}
&\phi_t = \frac{ (1+\phi_y^2) \phi_{xx} + (1+\phi_x^2) \phi_{yy} - 2\phi_x \phi_y \phi_{x,y}}{1 + \phi_x^2 + \phi_y^2}\\
&\phi(x,y,0) = g(x,y).
\end{split}
\right.
\end{equation}
As in two dimensions, we can obtain a Taylor expansion for the solution at time $t>0$:
\begin{equation}
\label{eq:exacttaylor3d}
\begin{split}
\phi(0,0,t)  &=t H(0,0) + \frac{1}{2} t^2 \left\{ \Delta^2 g(0,0) - 2H^3(0,0) + 6H(0,0) K(0,0) \right\} + O(t^3)\\
&= t H(0,0) + \frac{1}{2} t^2 \left\{  \Delta_S H(0,0) + H^3(0,0) - 2H(0,0) K(0,0) \right\} + O(t^3)
\end{split}
\end{equation}
where
\begin{equation}
\begin{split}
H(0,0) &= \Delta g(0,0) = g_{xx}(0,0) + g_{yy}(0,0)\\
K(0,0) &= g_{xx}(0,0) g_{yy}(0,0) - g^2_{xy}(0,0)\\
\Delta_S H(0,0) &= \Delta^2 g(0,0) - 3H^3(0,0) + 8 H(0,0) K(0,0).
\end{split}
\end{equation}
The Taylor expansion for the convolution step (\ref{eq:mbo1}) of Algorithm (\ref{eq:mbo1}) \& (\ref{eq:mbo2}) is now given by
\begin{equation}
\begin{split}
K * \mathbf{1}_{\Sigma^0}(0,0,z) =& \frac{1}{2} \theta(0) - \frac{1}{2\sqrt{\pi}} \frac{z}{\sqrt{t}} \theta(-1) + \frac{1}{24\sqrt{\pi}} \frac{z^3}{t^{3/2}} \theta(-3) + \frac{H(0,0)}{2\sqrt{\pi}} \sqrt{t} \theta(1)\\
&+ \frac{\Delta^2 g(0,0)}{4\sqrt{\pi}} t^{3/2} \theta(3) - \frac{1}{8\sqrt{\pi}} \frac{z^2}{\sqrt{t}} H(0,0) \theta(-1)\\
&+ \frac{1}{2\sqrt{\pi}} z \sqrt{t} \left( \frac{3}{4} H^2(0,0) - K(0,0) \right)  \theta(1)\\
&- \frac{1}{2\sqrt{\pi}} t^{3/2} \left( \frac{5}{4} H^3(0,0) - 3H(0,0)K(0,0) \right) \theta(3) + O(t^{5/2})
\end{split}
\end{equation}
Substituting the ansatz $z = a_1 t + a_2 t^2 + O(t^3)$, we get
\begin{equation}
K_t * \mathbf{1}_{\Sigma^0}(0,0,z) = A_0 + A_1 \sqrt{t} + A_3 t^{3/2} + O(t^{5/2})
\end{equation}
where
\begin{equation}
A_0 = \frac{\theta(0)}{2}
\end{equation}
and
\begin{equation}
A_1 = \frac{1}{2\pi} \Big( \theta(1) H(0,0) - a_1 \theta(-1) \Big)
\end{equation}
and
\begin{equation}
\begin{split}
A_3 =&- \frac{\theta(-1)}{2\sqrt{\pi}} a_2 + \frac{\theta(-3)}{24\sqrt{\pi}} a_1^3 + \frac{\theta(3)}{4\sqrt{\pi}} \Delta^2 g(0,0) - \frac{\theta(-1)}{8\sqrt{\pi}} a_1^2 H(0,0)\\
&+ \frac{\theta(1)}{2\sqrt{\pi}} a_1 \left( \frac{3}{4} H^2(0,0) - K(0,0) \right) - \frac{\theta(3)}{2\sqrt{\pi}} \left( \frac{5}{4} H^3(0,0) - 3H(0,0)K(0,0) \right)
\end{split}
\end{equation}
Choosing the thresholding level as 
\begin{equation}
\lambda = \frac{\theta(0)}{2} = \frac{1}{2} \sum_{j=1}^N c_j
\end{equation}
we set $A_1=0$ and solve for $a_1$ to obtain
\begin{equation}
a_1  = \frac{\theta(1)}{\theta(-1)} H(0,0)
\end{equation}
Having determined $a_1$, we substitute the expression for it into $A_3$ and solve for $a_2$ to obtain
\begin{equation}
\begin{split}
a_2 =& \frac{\theta(3)}{2\theta(-1)} \Delta^2 g(0,0)\\
&+ \frac{\theta^3(1) \theta(-3)+ 6 \theta^2(1) \theta^2(-1) - 15 \theta(3)\theta^3(-1)}{12 \theta^4(-1)} H^3(0,0)\\
&- \frac{\theta^2(1) - 3 \theta(3)\theta(-1)}{\theta^2(-1)}  H(0,0) K(0,0)
\end{split}
\end{equation}
Thus, in summary, the location of the interface along the $z$-axis after one time step with scheme (\ref{eq:mbo1}) \& (\ref{eq:mbo2}) is given by
\begin{equation}
z = t  B_1 H(0,0)  + t^2 \left\{ B_2 \Delta^2 g(0,0) + B_3 H^3(0,0) + B_4 H(0,0) K(0,0) \right\}
\end{equation}
where
\begin{equation}
\begin{split}
B_1 &= \frac{\theta(1)}{\theta(-1)},\\
B_2 &= \frac{\theta(3)}{2\theta(-1)}\\
B_3 &= + \frac{\theta^3(1) \theta(-3)+ 6 \theta^2(1) \theta^2(-1) - 15 \theta(3)\theta^3(-1)}{12 \theta^4(-1)} \mbox{, and}\\
B_4 &= - \frac{\theta^2(1) - 3 \theta(3)\theta(-1)}{\theta^2(-1)}
\end{split}
\end{equation}
To match the exact expansion (\ref{eq:exacttaylor3d}) for some possibly rescaled effective time step size, we need, in particular:
\begin{equation}
6 B_2 = B_4
\end{equation}
which, under the proviso that $\theta(-1) \not= 0$, implies $B_1 =0$.
That precludes matching (\ref{eq:exacttaylor3d}) up to $O(t^3)$ terms.
Hence, second order consistency with motion by mean curvature cannot be obtained using any linear combination of Gaussians as the convolution kernel in dimensions $d\geq 3$, even at the expense of violating the comparison principle (i.e. allowing the kernel to become negative).

\section{Conclusion}
We have exhibited a second order accurate in time scheme for curvature motion in the plane that is monotone: It preserves the comparison principle satisfied by the exact evolution it approximates.
The scheme is a variant of the threshold dynamics algorithm of Merriman, Bence, and Osher.
In particular, we have shown that there is a very special convolution kernel -- a carefully chosen linear combination of Gaussians --  to use in that algorithm that endows the scheme with both second order accuracy in time and monotonicity.
Numerical experiments presented bear out the advertised order of accuracy.
We have also shown that extending our work to three dimensions and higher will require a convolution kernel that cannot be as simple as a linear combination of Gaussians.
Some immediate, intriguing directions for further study include:
\begin{itemize}
\item Is there a more elaborate kernel construction that would result in a monotone, second order accurate algorithm of the form (\ref{eq:mbo1}) \& (\ref{eq:mbo2}) in three dimensions and higher?
\item The two dimensional special convolution kernel identified in Section \ref{sec:2d} is positive (which is what makes the resulting algorithm monotone), but its Fourier transform isn't.
This means the resulting scheme is not guaranteed to dissipate the Lyapunov functional identified in \cite{esedoglu_otto}.
Is there another kernel that results in second order accuracy in time, and that is positive in both physical and Fourier domains, so that both the comparison principle and energy based notions of stability are guaranteed?
\item We already know that the scheme is convergent, to the viscosity solution of curvature motion, thanks to consistency and monotonicity.
Given that consistency holds at second order,  can the rate of convergence be rigorously shown to be second order in time, as in \cite{kishii} that establishes first order convergence for the original algorithm?
\item Are there related geometric motions for which a similar thresholding scheme can be found that is monotone and second order?
\end{itemize}

\section{Acknowlegment}
Selim Esedo\=glu was supported by NSF DMS-2012015.
Jiajia Guo was supported by NSF DMS-1813003.

\bibliographystyle{plain}
\bibliography{main}

\begin{thebibliography}{10}

\bibitem{barles_georgelin}
G.~Barles and C.~Georgelin.
\newblock A simple proof of convergence for an approximation scheme for
  computing motions by mean curvature.
\newblock {\em SIAM J. Numer. Anal.}, 32:484--500, 1995.

\bibitem{barles_souganidis}
G.~Barles and P.~Souganidis.
\newblock Convergence of approximation schemes for fully nonlinear second order
  equations.
\newblock {\em Asymptotic Analysis}, 4:271--283, 1991.

\bibitem{bbc}
E.~Bonnetier, E.~Bretin, and A.~Chambolle.
\newblock Consistency result for a non-monotone scheme for anisotropic mean
  curvature flow.
\newblock {\em Interfaces and Free Boundaries}, 14(1):1--35, 2012.

\bibitem{chen_giga_goto}
Y.-G. Chen, Y.~Giga, and S.~Goto.
\newblock Uniqueness and existence of viscosity solutions of generalized mean
  curvature flow equations.
\newblock {\em Journal of Differential Geometry}, 33:749--786, 1991.

\bibitem{elsey_esedoglu_anisotropy}
M.~Elsey and S.~Esedo{\=g}lu.
\newblock Threshold dynamics for anisotropic surface energies.
\newblock {\em Mathematics of Computation}, 87(312):1721--1756, 2018.

\bibitem{ejz}
S.~Esedo\=glu, M.~Jacobs, and P.~Zhang.
\newblock Kernels with prescribed surface tension \& mobilty for threshold
  dynamics schemes.
\newblock {\em Journal of Computational Physics}, 337:62--83, 2017.

\bibitem{esedoglu_otto}
S.~Esedo{\=g}lu and F.~Otto.
\newblock Threshold dynamics for networks with arbitrary surface tensions.
\newblock {\em Communications on Pure and Applied Mathematics}, 68(5):808--864,
  2015.

\bibitem{esedoglu_ruuth_tsai2}
S.~Esedo{\=g}lu, S.~Ruuth, and Y.-H. Tsai.
\newblock Threshold dynamics for high order geometric motions.
\newblock {\em Interfaces and Free Boundaries}, 10(3):263--282, 2008.

\bibitem{evans_mbo}
L.~C. Evans.
\newblock Convergence of an algorithm for mean curvature motion.
\newblock {\em Indiana University Mathematics Journal}, 42:553--557, 1993.

\bibitem{evans_spruck}
L.~C. Evans and J.~Spruck.
\newblock Motion of level sets by mean curvature. {I}.
\newblock {\em Journal of Differential Geometry}, 33:635--681, 1991.

\bibitem{heintz05}
R.~Grzhibovskis and A.~Heintz.
\newblock A convolution-thresholding approximation of generalized curvature
  flows.
\newblock {\em SIAM Journal on Numerical Analysis}, 42(6):2652--2670, 2005.

\bibitem{grzibovskis_heintz}
R.~Grzhibovskis and A.~Heintz.
\newblock A convolution thresholding scheme for the {W}illmore flow.
\newblock {\em Interfaces and Free Boundaries}, 10(2):139--153, 2008.

\bibitem{ishii_pires_souganidis}
H.~Ishii, G.~E. Pires, and P.~E. Souganidis.
\newblock Threshold dynamics type approximation schemes for propagating fronts.
\newblock {\em Journal of the Mathematical Society of Japan}, 51:267--308,
  1999.

\bibitem{kishii}
K.~Ishii.
\newblock Optimal rate of convergence of the {B}ence-{M}erriman-{O}sher
  algorithm for motion by mean curvature.
\newblock {\em SIAM Journal on Numerical Analysis}, 37(3):841--866, 2005.

\bibitem{laux_otto}
T.~Laux and F.~Otto.
\newblock Convergence of the thresholding scheme for multi-phase mean-curvature
  flow.
\newblock {\em Calculus of Variations and Partial Differential Equations},
  55(5):1--74, 2016.

\bibitem{laux_otto2}
T.~Laux and F.~Otto.
\newblock Brakke's inequality for the thresholding scheme.
\newblock {\em Calculus of Variations and Partial Differential Equations},
  59(1):39--65, 2020.

\bibitem{mascarenhas}
P.~Mascarenhas.
\newblock Diffusion generated motion by mean curvature.
\newblock CAM Report 92-33, UCLA, July 1992.
\newblock (URL = http://www.math.ucla.edu/applied/cam/index.html).

\bibitem{merriman_bence_osher}
B.~Merriman, J.~K. Bence, and S.~Osher.
\newblock Motion of multiple junctions: a level set approach.
\newblock {\em Journal of Computational Physics}, 112(2):334--363, 1994.

\bibitem{mbo92}
B.~Merriman, J.~K. Bence, and S.~J. Osher.
\newblock Diffusion generated motion by mean curvature.
\newblock In J.~Taylor, editor, {\em Proceedings of the Computational Crystal
  Growers Workshop}, pages 73--83. AMS, 1992.

\bibitem{ruuththesis}
S.~Ruuth.
\newblock {\em Efficient algorithms for diffusion-generated motion by mean
  curvature}.
\newblock PhD thesis, University of British Columbia, 1996.

\bibitem{ruuth0}
S.~J. Ruuth.
\newblock Efficient algorithms for diffusion-generated motion by mean
  curvature.
\newblock {\em Journal of Computational Physics}, 144:603--625, 1998.

\bibitem{ruuth_merriman}
S.~J. Ruuth and B.~Merriman.
\newblock Convolution generated motion and generalized {H}uygens' principles
  for interface motion.
\newblock {\em SIAM Journal on Applied Mathematics}, 60:868--890, 2000.

\bibitem{swartz}
D.~Swartz and N.~K. Yip.
\newblock Convergence of diffusion generated motion to motion by mean
  curvature.
\newblock {\em Communications in Partial Differential Equations},
  42(10):1598--1643, 2017.

\bibitem{walkington}
N.~Walkington.
\newblock Algorithms for computing motion by mean curvature.
\newblock {\em SIAM Journal on Numerical Analysis}, 33(6):2215--2238, 1996.

\bibitem{zaitzeff1}
A.~Zaitzeff, S.~Esedo\=glu, and K.~Garikipati.
\newblock Second order threshold dynamics schemes for two phase motion by mean
  curvature.
\newblock {\em Journal of Computational Physics}, 410, 2020.

\end{thebibliography}

\end{document}